\documentclass[12pt, reqno]{amsart}
\usepackage{amsmath,amssymb,amsthm}
\usepackage{color}
\usepackage[colorlinks, linkcolor=red, citecolor=blue, urlcolor=blue, hypertexnames=false]{hyperref}

\newtheorem{lemma}{Lemma}[section]

\newtheorem{theorem}{Theorem}[section]

\theoremstyle{definition}

\theoremstyle{remark}

\theoremstyle{remark}
\newtheorem{remark}{Remark}[section]
\numberwithin{equation}{section}

\newcommand{\R}{{\mathbb R}}

\definecolor{blu}{rgb}{0,0,1}

\title[Existence and Asymptotic Behavior of Solutions on Graphs]{Existence and Asymptotic Behavior of Solutions to Nonlinear Bi-harmonic Equations on Locally Finite Graphs}

\author[Tingjian Luo]{Tingjian Luo*}
\address{Tingjian Luo
\newline\indent
School of Mathematics and Information Sciences, 
\newline\indent
Guangzhou University, Guangzhou, China}
\email{luotj@gzhu.edu.cn}

\author[Yongqiang Wei]{Yongqiang Wei}
\address{Yongqiang Wei
\newline\indent
School of Mathematics and Information Sciences, 
\newline\indent
Guangzhou University, Guangzhou, China}
\email{}

\begin{document}


\begin{abstract} In this paper, we investigate the existence and asymptotic behaviors of non-trivial solutions for a class of nonlinear bi-harmonic equations on locally finite graphs \(G=(V,E)\):
\begin{equation}\label{eq1.0}
\Delta^2 u-\Delta u+(a(x)+1)u=f(u)+\varepsilon h(x).
\end{equation}
where the potential \(a(x): V\to R\) is bounded from below but may change sign. Using mainly variational methods, we prove the multiplicity of non-trivial solutions of the equations on locally finite graphs. Among them, one solution is of local minimizer type with negative energy and another is of mountain pass type with positive energy. In addition, we demonstrate the asymptotic behaviors of the solutions we obtained as the parameter \(\varepsilon\to 0^+\). 
\end{abstract}
\thanks{*Corresponding author: luotj@gzhu.edu.cn (Tingjian Luo)}
\maketitle

\section{Introduction}
Let \(G=(V,E)\) be a discrete graph, where \(V\) denotes the set of vertices and \(E\) the set of edges. For any \(x\in V\), there are only finite \(y\in V\) such that \(xy\in E\), then \(G\) is called a locally finite graph. \(G\) is called connected if any two vertices \(x\) and \(y\) can be connected via finite edges. If the set of vertices \(V\) is finite, then \(G=(V,E)\) is called a bounded domain, and otherwise an unbounded domain. In this paper, we assume that \(G\) is a locally finite graph and is connected (which may contain infinite vertices, becoming unbounded domain). For a measure \(\mu: V\to R^{+}\) defined on \(V\), we assume that there exists a \(\mu_{min}>0\), such that \(\mu(x)\geq \mu_{min}\) for all \(x\in V\). Let \(w_{xy}\) be the weight of an edge \(xy\in E\), such that \(w_{xy}>0\) and \(w_{xy}=w_{yx}\).

Graphs, as natural structures, are widely applied in modeling problems across diverse fields such as image processing, neural networks \cite{BL}, and social networks. In recent years, research on partial differential equations on graphs has evolved into a highly active area, with particular focus on equations in geometric analysis and mathematical physics. This has become a crucial foundation for advancing related analytical work. Based on the definition on the graph, many authors have studied the existence of solutions. For instance, in \cite{LW} the authors investigated the existence and non-existence of global solutions for a semi-linear heat equation on graphs. In \cite{YLY2}, Grigo'yan et al. investigated the existence of positive solutions for certain nonlinear equations on locally finite graphs. Also there are some work on the Kazdan-Warner equations on graphs, see e.g. \cite{CM, G, YLY1, SW} .

The bi-harmonic nonlinear Schr\"odinger equations defined on the Euclidian space has been well studied in the past decades, which appeared  as a canonical model for nonlinear Hamiltonian PDEs with dispersion of super-quadratic order. Historically the first work should go back to Karpman and Karpman-Shagalov in the physics literature \cite{K, KS}. Since then, more and more interest has been attracted around the analyses related to the existence, non-existence, stability, instability and blow-up of solutions to the fourth-order equations with a second order dispersion. See for example \cite{BCGJ1,BCGJ2,BCM,BCSN,BL,BN, FJMM, GHZ, S, WS, ZLJS, ZZ} , also more recently \cite{CF} and the references therein.   

Concerned with the study on the Eq.\eqref{eq1.0} on the locally finite graphs, we may refer to the work \cite{HSZ}, where the authors established the existence of ground state solutions of Eq. \eqref{eq1.0} with \(\varepsilon=0\) on an infinite graph, and investigated some concentration toward a bounded potential well for the solutions. As the case \(\varepsilon > 0\), the authors in \cite{H} proved on a bounded domain the existence of two solutions for a class of second order equation as the parameter \(\varepsilon>0\) small enough. In the case of bounded domain, one may enjoy some benefit of compact embedding in establishing the compactness of approximating solutions. In the unbounded case, with the assumption of positive coercive potential, the authors established in \cite{YLY2} the existence of two positive solutions. We also note a recent work \cite{LZ} where the authors established a multiple result of solutions for a class of second order equations with a perturbation and a coercive and positive potential.

The aim of this paper is twofold. First, we prove the existence of two solutions of Eq. \eqref{eq1.0} as \(\varepsilon>0\) small enough, where one is of local minimizer type with negative energy and the other is of mountain-pass type with positive energy. This will extend the work in \cite{H} to the case of unbounded domain, and also we do not need to assume as in \cite{LZ,YLY2} the potential to be coercive and positive, which is even allowed to be sign-changing. Second, we establish the asymptotic behavior of the solutions obtained as the parameter \(\varepsilon \to 0^+\).

To begin with, let us introduce some notations on the analysis on locally finite graphs. For any function \( u(x): V \to \R \) and \( x \in V \), the \( \mu \)-Laplacian (or Laplacian for short) of \( u \) at \( x \) is defined by
\[
\Delta u(x) := \frac{1}{\mu(x)} \sum_{y \sim x} \omega_{xy} \bigl( u(y) - u(x) \bigr).
\]
The gradient form \( \Gamma(u, v) \) of two functions \( u \) and \( v \) at \( x \in V \) is given by
\[
\Gamma(u, v)(x) := \frac{1}{2\mu(x)} \sum_{y \sim x} \omega_{xy} \bigl( u(y) - u(x) \bigr) \bigl( v(y) - v(x) \bigr).
\]
For brevity, we use \( \Gamma(u) \) instead of \( \Gamma(u, u) \) and sometimes \( \nabla u \nabla v \) instead of \( \Gamma(u, v) \). The length of \( \Gamma(u) \) at \( x \in V \) is denoted by
\[
|\nabla u(x)|:= \sqrt{\Gamma(u)(x)} = \left( \frac{1}{2\mu(x)} \sum_{y \sim x} \omega_{xy} \bigl( u(y) - u(x) \bigr)^2 \right)^{\frac{1}{2}}.
\]
For any function \( u: V \to \R \), an integral of \( u \) over \( V \) is defined by
\[
\int_V u d\mu := \sum_{x \in V} \mu(x) u(x).
\]
The bi-harmonic operator of \(u:V\to\R\), denoted by \(\Delta^{2}u\), is defined in the distributional sense as
\[
\int\limits_{V}(\Delta^{2}u)\phi\,d\mu:=\int\limits_{V}\Delta u\Delta\phi\,d\mu,
\quad\forall\phi\in C_{c}(V),
\]
where \(C_{c}(V):=\Big\{u:V\to \R: \mbox{ the set } \{x\in V:u(x)\neq 0\} \mbox{ has finite cardinality } \Big\}\). \\

Let \( W^{2,2}(V) \) be the completion of \( C_c(V) \) under the norm
\[
\|u\|_{W^{2,2}(V)} := \left( \int_V \bigl( |\Delta u|^2 + |\nabla u|^2 + u^2 \bigr) d\mu \right)^{\frac{1}{2}}.
\]
Clearly, \( W^{2,2}(V) \) is a Hilbert space with the inner product
\[
(u, v)_{W^{2,2}(V)} := \int_V \bigl( \Delta u \Delta v + \nabla u \nabla v + uv \bigr) d\mu, \ \forall u, v \in W^{2,2}(V).
\]
It is natural to consider the weighed space
\[
H:= \left\{ u \in W^{2,2}(V) : \int_V a(x)u^2 d\mu < +\infty \right\},
\]
endowed with the norm
\[
\|u\|_{H} := \left( \int_V \bigl( |\Delta u|^2 + (a(x) + 1) u^2 \bigr) d\mu \right)^{\frac{1}{2}}.
\]
The space \(H\) is also a Hilbert space with its inner product
\[
(u,v)_{H} := \int_V \bigl( \Delta u \Delta v+ (a(x) + 1) uv \bigr) d\mu, \ \forall u, v \in H.
\]

In this paper, we are concerned with the following bi-harmonic nonlinear Sch\"odinger equation with a potential and an \(\varepsilon-\)perturbation:
\begin{equation}\label{eq1.1}
\Delta^2 u-\Delta u+(a(x)+1)u=f(u)+\varepsilon h(x),\ x\in V,  
\end{equation}
where the assumptions on the functions $a(x), f(s)$ and $h(x)$ will be given in the sequel. 

Let $\varepsilon=0$, then Eq.\eqref{eq1.1} is reduced to its limit form:
\begin{equation}\label{eq1.2}
\Delta^2 u-\Delta u+(a(x)+1)u=f(u),\ x\in V.   
\end{equation}

The corresponding energy functional \(J_\varepsilon: H\to \R\) related with Eq.\eqref{eq1.1} is 
\begin{equation*}
J_\varepsilon(u) :=\frac{1}{2}\int_{V}\Big(|\Delta u|^2+|\nabla u|^2+(a(x)+1)|u|^2\Big)d\mu-\int_{V}F(u)d\mu-\int_{V}\varepsilon h(x)ud\mu.
\end{equation*}
Then, from the point view of variational methods, a weak solution \(u\) to Eq.\eqref{eq1.2} corresponds to a critical point of \(J_\varepsilon\) on \(H\). Namely, for any test function \(\varphi\in H\),
\begin{equation*}
(J'_\varepsilon(u),\varphi)=\int_{V}[\Delta u\Delta \varphi+\nabla u\nabla \varphi+(a(x)+1)u\varphi-f(u)\varphi-\varepsilon h(x)\varphi]d\mu=0.
\end{equation*}

Therefore in this paper we will assume always that\\
\((A_1)\) The potential \(a(x): V\to \R\), satisfies \(\inf\limits_{x\in V}a(x)\geq a_0\) for some \(a_0\in(-1,0)\);\\
\((A_2)\) For every \(M>0\),  the volume of the set \(D_M\) is finite, namely,
\begin{equation*}
Vol(D_M)=\sum\limits_{x\in D_M}\mu(x)<\infty,
\end{equation*}
where \(D_M:=\{x\in V: a(x)\leq M\}\);\\
\((A'_2)\) \((a(x)+1)^{-1}\in L^1(V)\);\\
\((f_1)\) \(f(s)\in C(\R,\R)\), and \(\lim\limits_{s\to 0}\frac{f(s)}{s}=0\);\\
\((f_2)\) \(\lim\limits_{|s|\to +\infty}\frac{f(s)}{|s|^p}=0\), as \(1<p<+\infty\);\\
\((f_3)\) There is constant \(\theta >2\), such that
\begin{equation}\label{1.4}
 0<\theta F(s)\leq sf(s),\qquad\forall\ |s|>0,   
\end{equation}
where \(F(s):=\int_{0}^{t}f(t)dt\);\\
\((h_1)\) \(h(x)\in L^2(V)\) and \(\|h\|_2 \neq 0\).\\ 

Concerned with the above assumptions on the functions \(a(x), f(s), h(x)\), we would like to give some remarks:
\begin{remark}
\begin{itemize}
  \item [(1)] From \((f_1)\) and \((f_2)\), one may deduce that for any \(\alpha\), there exists constant \(C(\alpha)>0\) such that
\begin{equation}\label{1.3}
|f(s)|\leq \alpha |s|+C(\alpha)|s|^p,\quad \forall\ |s|>0.   
\end{equation}
The assumption \eqref{1.4} is referred as the \emph{Ambrosetti–Rabinowitz} condition, from which we derive that
\begin{equation}\label{1.5}
\quad F(s)\geq \beta\cdot |s|^{\theta},\quad \forall\ |s|>1.  
\end{equation}
where \(\beta>0\) independent of \(s\), see \cite{AR}.
  \item [(2)] Under the assumption \((A_2)\) or \((A'_2)\), we are able to establish the compact embedding of the space \(H\) to \(L^q(V)\), see Lemma \ref{lm2.1} and Lemma \ref{2.2}, which extend the compact embedding theorems in \cite {CWY} to the bi-harmonic setting. Note particularly that, in the Euclidian space, the condition \((A_2)\) is \textbf{weaker} than the coercivity condition (namely \(a(x)\to +\infty\) as \(|x|\to +\infty\)), however, in the locally finite graph, it is easy to check that they are equivalent to each other under the assumption that \(\mu(x)\geq\mu_{min}>0\) for all \(x\in V\).
  \item [(3)] Both functions \(a(x)\) and \(h(s)\) are allowed to be \textbf{sign-changing}.
\end{itemize}

\end{remark}

Now, we state our first main result.

\begin{theorem}\label{th1.1}
Assume that \((A_1)\), \((A_2)\) or \((A'_2)\), \((f_1)-(f_3)\) and \((h_1)\) hold. Then as the parameter \(\varepsilon>0\) sufficiently small,  Eq.\eqref{eq1.1}  have two non-trivial solutions \(u_{\varepsilon}^{(1)}\) and \(u_{\varepsilon}^{(2)}\), where \(u_{\varepsilon}^{(1)}\) is of mountain pass type with positive energy and \(u_{\varepsilon}^{(2)}\) is of local minimizer type with negative energy.
\end{theorem}

\begin{remark}
Since the graph we consider is possibly infinite, we establish a new compact embedding result based on the potential to recover the compactness of approximating solutions, see Lemma \ref{lm2.1} and Lemma \ref{lm2.2}, but where the potential can be sign-changing. The challenging part is then to establish the mountain-pass geometry and local minimizing structure, as the parameter \(\varepsilon>0\) sufficiently small, see Lemma \ref{lm3.3.3} and Lemma \ref{lm3.4} respectively. 
\end{remark}

\begin{remark}
We also point out that in \cite{YLY2} the authors proved further that the solutions obtained are strictly positive in the case of second order regime. Their strategy is using firstly the negative part \(u^-\) of solutions as a testing function to prove that \(u^-(x)=0, \forall x\in V\), then arguing by contradiction to prove that the solution \(u(x)\neq 0, \forall x\in V\), and finally concluding that \(u(x)> 0, \forall x\in V\). This strategy is not applicable for the bi-harmonic equations. One of the reasons is that \(u^-\) is possibly not in the space \(H\), and also it is hard to determine the sign of term related to the fourth-order operator. Hence, to know the exact sign of the solutions obtained in this theorem, we may need to develop a new approach.    
\end{remark}

Our second existence result is concerned with the case \(\varepsilon=0\). We prove the existence of a mountain-pass solution to the limit equation Eq.\eqref{eq1.2}.

\begin{theorem}\label{th1.2}
Assume that \((A_1)\), \((A_2)\) or \((A'_2)\), \((f_1)-(f_3)\) hold. Then Eq.\eqref{eq1.2} has at least one non-trivial solution \(u_0\in H\) with its mountain-pass level \(c_0>0\) given by \eqref{mp_0}.
\end{theorem}

Similar to Theorem \ref{th1.1}, we define the energy functional \(J_0:H\to \R\) with respect to \eqref{eq1.2}
\begin{equation*}
J_0(u) :=\frac{1}{2}\int_{V}\Big(|\Delta u|^2+|\nabla u|^2+(a(x)+1)u^2\Big)d\mu-\int_{V}F(u)d\mu.
\end{equation*}

Finally, we investigate the asymptotic behavior of the two solutions obtained in Theorem \ref{th1.1} as the parameter \(\varepsilon\to 0^+\) and the result is the following: 

\begin{theorem}\label{th1.3}
Under the conditions given in Theorem \ref{th1.1} and Theorem \ref{th1.2}, if \(\varepsilon\to 0^+\), then 
\[u_{\varepsilon}^{(1)} \to u_1, \ \mbox{in}\ H,\  \mbox{and}\ u_{\varepsilon}^{(2)} \to 0, \ \mbox{in}\ H. \]
where \(u_1\in H\) is a mountain pass solution of Eq.\eqref{eq1.2} at the energy level \(J_0(u_1)=c_0\). In particular, 
\[c_\varepsilon=J_{\varepsilon}(u_{\varepsilon}^{(1)}) \to J_0(u_1)=c_0,\] 
where \(c_{\varepsilon}, c_0\) are given by \eqref{mp} and \eqref{mp_0} respectively. 
\end{theorem}

\begin{remark}
The point that \(c_{\varepsilon} \to c_0\) as \(\varepsilon\to 0^+\), is not trivial in generally. Since the mountain pass \(\|\gamma(t)\|_2\) for \(\gamma\in \mathcal{L}\) is not necessarily uniformly bounded. For example, one can find a sequence of \(\gamma_n(t):=nt(1-t)+tu'\) (where \(u'\in H\) is given in Lemma \ref{lm3.2}), then \(\gamma_n\in \mathcal{L}\) but clearly \(\|\gamma_n(t)\|_2\) is not bounded. To overcome this difficulty, we define the two mountain pass levels \(c_{\varepsilon}\) and \(c_0\) under the same path set \(\mathcal{L}\), and develop some careful analyses on the energy estimates, see Lemma \ref{lm4.2} and its proof for more details.
\end{remark}

\section{Some Preliminary Results}
To prove Theorem \ref{th1.1} and Theorem \ref{th1.2}, we need to employ the Sobolev embedding theorem on locally finite graphs. This leads us to introduce the following lemma:
\begin{lemma}\label{lm2.1}
Assume that the conditions \((A_1)\) and \((A_2)\) hold. Then the space \(H\) is weakly compact, and \(H\) is compactly embedded in \(L^q(V)\), for all \(2\leq q\leq+\infty\).
\end{lemma}

\begin{proof}
For any \(x_0\in V\) and \(u\in H\), under the condition \((A_1)\), we have:
\begin{equation*}
\begin{aligned}
\| u \|_{H}^2&=\int_{V}\Big(|\Delta u|^2+(a(x)+1)|u|^2\Big)d\mu\\
&\geq\int_{V}(a(x)+1)u^2d\mu\\
&\geq(a_0+1)\sum\limits_{x\in V}\mu(x)|u(x)|^2\\
&\geq(a_0+1)\mu_{min}|u(x_0)|^2,
\end{aligned}
\end{equation*}
thus,
\begin{equation}\label{2.1}
|u(x_0)|\leq\sqrt{\frac{1}{(a_0+1)\mu_{min}}}\| u \|_{H}.   
\end{equation}
Therefore, \(H\) is continuously embedded into \(L^\infty(V)\). According to the interpolation inequality, there exists constant $C_{\infty}>0$, such that 
\begin{equation*}
\| u \|_{L^q}\leq\| u \|_{L^2}^\chi\| u \|_{L^\infty}^{1-\chi}\leq C_{\infty}\| u \|_{H},\quad \forall u\in H,
\end{equation*}
for all \(2\leq q<\infty\), and where $\chi\in (0,1)$. Thus, we conclude that \(H\) is continuously embedded into \(L^q(V)\) for all \(2\leq q\leq\infty\). 

Now, we prove that \(H\) is compactly embedded in \(L^q(V)\), for all \(2\leq q\leq+\infty\). To this aim, let  \(\{u_k\}\) be a bounded sequence in \(H\), since \(H\) is a Hilbert space and reflexive (see \cite{HSZ}), then there exist a subsequence (still denoted by \(\{u_k\}\)) and some \(u\in H\) such that \(u_k\rightharpoonup u\) in \(H\). In particular, one has that \(u_k\rightharpoonup u\) in \(L^2(V)\). Then, for any \(\varphi\in L^2(V)\),
\begin{equation}\label{2.2}
\lim\limits_{k \to \infty}\int_{V}(u_k-u)\varphi d\mu=\lim\limits_{k \to \infty}\sum\limits_{x\in V}\mu(x)(u_k(x)-u(x))\varphi(x)=0.   
\end{equation}
For any \(x_0\in V\), let
\begin{equation*}
\varphi_0(x) := \left\{
\begin{array}{ll}
1, & \text{if } x = x_0;\\
0, & \text{if } x \neq x_0.
\end{array}
\right.
\end{equation*}
Clearly \(\varphi_0(x)\in L^2(V)\), substituting \(\varphi_0\) into \eqref{2.2} yields \(\lim\limits_{k \to \infty}u_k(x_0)=u(x_0)\). Since \(x_0\) is arbitrary, then \(\lim\limits_{k \to \infty}u_k(x)=u(x)\) for any \(x \in V\).

Now we prove that \(u_k\to u\) in \(L^p(V)\) for all \(2\leq p\leq \infty\). Indeed, by the interpolation inequality and the following estimate
\[\|u_k-u\|_{L^{\infty}}^2\leq \frac{1}{\mu_{min}}\int_V|u_k-u|^2d\mu,\]
it is enough to verify that \(u_k\to u\) is in \(L^2(V)\). Namely,
\begin{equation}\label{2.3}
\lim\limits_{k \to \infty}\int_{V}|u_k-u|^2d\mu=0.  
\end{equation}
Let \(D_M\) be given in \((A_2)\), then 
\begin{equation*}
\int_{V}|u_k-u|^2d\mu=\int_{D_M}|u_k-u|^2d\mu+\int_{V\backslash D_M}|u_k-u|^2d\mu.
\end{equation*}
On one hand, since \(\{\|u_k-u\|_{H}\}\) is bounded, then as $M\to +\infty$, we have
\begin{equation*}
\int_{V\backslash D_M}|u_k-u|^2d\mu\leq\int_{V\backslash D_M}\frac{a(x)}{M}|u_k-u|^2d\mu\leq\frac{\|u_k-u\|_{H}^2}{M}\to 0,
\end{equation*}
On the other hand, let \(V_1 : =\{{x\in V : d(x,x_0)\leq R}\}\), and
\begin{equation*}
\int_{D_M}|u_k-u|^2d\mu=\int_{D_M\cap V_1}|u_k-u|^2d\mu+\int_{D_M\backslash V_1}|u_k-u|^2d\mu.
\end{equation*}
By the H\"older inequality, we have
\begin{equation*}
\begin{aligned}
\int_{D_M\backslash V_1}|u_k-u|^2d\mu&\leq(\int_{D_M\backslash V_1}|u_k-u|^{2p}d\mu) ^{\frac{1}{p}}(\int_{D_M\backslash V_1}1d\mu)^{\frac{1}{q}}\\&\leq\|u_k-u \|_{2p}^2(\sum\limits_{x\in D_M\backslash V_1}\mu(x))^{\frac{1}{q}},
\end{aligned}
\end{equation*}
where \(\frac{1}{p}+\frac{1}{q}=1\). In view of \(Vol (D_M)<\infty\), then as \(R\) large enough, we have
\begin{equation*}
Vol(D_M\backslash V_1)=\sum\limits_{x\in D_M\backslash V_1}\mu(x)\leq(\frac{1}{R})^q,
\end{equation*}
Since \(\{\|u_k-u\|_{2p}\}\) is bounded, then letting \(R\to +\infty\),
\begin{equation*}
\int_{D_M\backslash V_1}|u_k-u|^2d\mu\leq\|u_k-u\|_{2p}^2(\sum\limits_{x\in D_M\backslash V_1}\mu(x))^{\frac{1}{q}}\leq \frac{\|u_k-u\|_{2p}^2}{R}\to 0,
\end{equation*}
By definition, \(V_1\) is a bounded domain, thus
\begin{equation*}
\int_{D_M\cap V_1}|u_k-u|^2d\mu\to 0.
\end{equation*}
Therefore, we conclude that
\begin{equation*}
\int_{D_M}|u_k-u|^2d\mu\to 0,
\end{equation*}
Thus, \eqref{2.3} follows, and the proof is completed.
\end{proof}

\begin{remark}\label{rek2.1}
By Lemma \ref{lm2.1}, for any \(2\leq q\leq +\infty\), there exists a constant \(C_q\) depending only on \(q\) such that
\begin{equation}\label{2.4}
\|u\|_{q}\leq C_q\|u\|_H.   
\end{equation}
\end{remark}

Replaced \((A_2)\) by \((A'_2)\), we have the following compact embedding result.

\begin{lemma}\label{lm2.2}
Assume that the conditions \((A_1)\) and \((A'_2)\) holds. Then the space \(H\) is weakly compact, and \(H\) is compactly embedded in \(L^q(V)\), holding for all \(1\leq q\leq+\infty\).
\end{lemma}

\begin{proof}
As the same as the proof in Lemma \ref{lm2.1}, under the condition \((A_1)\), we can prove that \(H\) is continuously embedded into \(L^\infty(V)\). Next, we prove that when \(1\leq q<\infty\), \(H\hookrightarrow L^q(V)\). From \((A'_2)\), for any \(u\in H\), we have
\begin{equation*}
\begin{aligned}
\int_{V}|u|d\mu&=\int_{V}(a(x)+1)^{-\frac{1}{2}}(a(x)+1)^{\frac{1}{2}}|u|d\mu\\
&\leq\Big(\int_{V}(a(x)+1)^{-1}d\mu\Big)^{\frac{1}{2}}\Big(\int_{V}(a(x)+1)|u|^2d\mu\Big)^{\frac{1}{2}}\\
&=\|(a(x)+1)^{-1}\|_{L^1}^{\frac{1}{2}}\Big(\int_{V}(a(x)+1)|u|^2d\mu\Big)^{\frac{1}{2}}\\
&\leq\|(a(x)+1)^{-1}\|_{L^1}^{\frac{1}{2}}\|u\|_H.
\end{aligned}
\end{equation*}
This verifies that \(H\hookrightarrow L^1(V)\), and thus by the interpolation inequality, we finally conclude that \(H\) can be continuously embedded into \(L^q(V)\) for all \(1\leq q\leq\infty\).  

To prove that the above embedding is compact, as the proof in Lemma \ref{lm2.1}, we only need to verify that any bounded sequence \(\{u_k\}\) in \(H\) has a weakly convergent subsequence in \(L^1(V)\). Up to a subsequence, there exists \(u\) such that \(u_k\rightharpoonup u\) in \(H\). In particular, \(u_k\rightharpoonup u\) in \( L^1(V)\).


Now, we prove that \(u_k\to u\) in \(L^1(V)\). Since \(\{u_k\}\) is bounded in \(H\), there exists a constant \(C>0\) independent of $k$, such that
\begin{equation*}
\int_{V}(a(x)+1)|u_k-u|^2d\mu\leq C.
\end{equation*}
For fixed \(x_0\in V\), for any \(\varepsilon>0\), it follows from \((A'_2)\) that there exists \(R>0\) such that
\begin{equation*}
\int_{d(x,x_0)>R}(a(x)+1)^{-\frac{1}{2}}d\mu<\varepsilon^2.
\end{equation*}
Then, using the H\"older inequality, we have
\begin{equation*}
\begin{aligned}
\int_{d(x,x_0)>R}|u_k-u|d\mu&=\int_{d(x,x_0)>R}(a(x)+1)^{-\frac{1}{2}}(a(x)+1)^{\frac{1}{2}}|u_k-u|d\mu\\
&\leq\Big(\int_{d(x,x_0)>R}(a(x)+1)^{-1}d\mu\Big)^{\frac{1}{2}}\Big(\int_{d(x,x_0)>R}(a(x)+1)|u_k-u|^2d\mu\Big)^{\frac{1}{2}}\\&\leq C^{\frac{1}{2}}\varepsilon.
\end{aligned}
\end{equation*}
Since \(\varepsilon>0\) is arbitrary, the above estimate yields that 
\begin{equation*}
\begin{aligned}
\lim\limits_{k\to +\infty}\int_{d(x,x_0)>R}|u_k-u|d\mu=0.
\end{aligned}
\end{equation*}
In view of 
\begin{equation*}
\lim\limits_{k\to +\infty}\int_{d(x,x_0)\leq R}|u_k-u|d\mu=0,
\end{equation*}
then we conclude that
\begin{equation*}
\liminf\limits_{k\to +\infty}\int_{V}|u_k-u|d\mu=0.
\end{equation*}
Namely, \(u_k\to u\) in \(L^1(V)\). 
At this point, the proof is completed.
\end{proof}


\section{Proofs of Theorem \ref{th1.1} and Theorem \ref{th1.2}}

In this section, we give the proofs of Theorem \ref{th1.1} and Theorem \ref{th1.2}, concerned with the existence results. First, we prove the mountain pass geometry of the functional \(J_\varepsilon(u)\) on \(H\). 
\begin{lemma}\label{lm3.1}
For given \(\varepsilon\in (0,1)\) small enough, there exist two constants \(0<r_\varepsilon<1\) and \(M_\varepsilon>0\), such that \(J_\varepsilon(u)\geq M_\varepsilon\), for any \( u\in H \) satisfying \(\| u \|_{H}=\frac{1}{2}r_\varepsilon\) or \( r_\varepsilon\).
\end{lemma}

\begin{proof}
Using \eqref{1.3}, \eqref{2.4}, and the H\"older inequality, we have
\begin{align}\label{estima}
\begin{split}
J_\varepsilon(u) &=\frac{1}{2}\int_{V}\Big(|\Delta u|^2+|\nabla u|^2+(a(x)+1)|u|^2\Big)d\mu-\int_{V}F(u)d\mu-\int_{V}\varepsilon h(x)ud\mu \\ 
&=\frac{1}{2}\| u \|_{H}^2-\int_{V}F(u)d\mu-\int_{V}\varepsilon h(x)ud\mu+\frac{1}{2}\|\nabla u\|_{2}^2\\
& \geq\frac{1}{2}\| u \|_{H}^2-\frac{\alpha}{2}\int_{V}|u|^2d\mu-\frac{C(\alpha)}{p+1}\int_{V}| u |^{p+1}d\mu-\varepsilon\|h\|_2 \|u\|_2\\
& \geq(\frac{1}{2}-\frac{\alpha}{2(a_0+1)})\| u \|_{H}^2-\frac{C_{p+1}}{p+1}C(\alpha)\| u \|_{H}^{p+1}-\varepsilon (a_0+1)^{-\frac{1}{2}}\|h\|_2 \|u\|_{H}.
\end{split}
\end{align}
Let \(\alpha=\frac{a_0+1}{2}\) and \(r_\varepsilon=\sqrt{\varepsilon}<1\), then for any \( u\in H \) satisfying \(\| u \|_{H}=\delta r_\varepsilon \) with \(\delta:=\frac{1}{2} \text{ or } 1\), we have
\begin{align*}
J_\varepsilon(u)&\geq \frac{\delta^2}{4}\varepsilon-\delta^{p+1}\frac{C_{p+1}}{p+1}C(\alpha)\varepsilon^{\frac{p+1}{2}}-\delta\varepsilon^{\frac{3}{2}} (a_0+1)^{-\frac{1}{2}}\|h(x)\|_{2}\\
&\geq \frac{1}{16}\varepsilon-\frac{C_{p+1}}{p+1}C(\alpha)\varepsilon^{\frac{p+1}{2}}-\varepsilon^{\frac{3}{2}} (a_0+1)^{-\frac{1}{2}}\|h(x)\|_{2}.
\end{align*}
Set \(M_\varepsilon:=\frac{1}{32}\varepsilon\) and let \(\varepsilon>0\) be sufficiently small, then we have \(J_\varepsilon(u)\geq M_\varepsilon >0\).
\end{proof}

Next, we prove a technique lemma used in the sequel.

\begin{lemma}\label{lm-tech}
Assume that \(h(x)\) satisfies \((h_1)\), then there exists \(\hat{u} \in H\), such that 
\[\int_V h(x)\hat{u}(x)d\mu>0.\]
\end{lemma}
\begin{proof}
The idea of the proof partially stems from \cite[Lemma 2.3]{H}. Let us consider the following auxiliary equation:
\begin{equation}\label{eq3.1}
\Delta^2 u+[a(x)+1]u= h(x), x\in V, 
\end{equation}
whose corresponding energy functional is given by
\begin{equation*}
E_{h}(u):=\frac{1}{2}\| u \|_{H}^2-\int_{V}h(x)u(x)d\mu, \quad \forall u\in H.
\end{equation*}
Then, as \eqref{estima}, we have
\begin{equation}\label{3.3}
E_{h}(u)\geq\frac{1}{2}\| u \|_{H}^2- (a_0+1)^{-\frac{1}{2}}\|h(x)\|_{2}\|u\|_{H},\quad \forall u\in H.
\end{equation}
This shows that the functional \(E_h(u)\) is bounded from below on \(H\). Thus, we can define
\begin{equation*}
N_h:=\inf\limits_{u\in H}E_h(u).
\end{equation*}
Let \(\{u_k\}\in H\) be a minimizing sequence of \(N_h\), namely \(u_k\in H\) satisfies \(\lim\limits_{k\to \infty} E_{h}(u_k)= N_h\). 
From \eqref{3.3}, we observe that \(\{u_k\}\) is bounded in \(H\). Then there exists \(\hat{u} \in H\) such that \(u_k\rightharpoonup \hat{u}\) in \(H\). Furthermore, by Lemma \ref{lm2.1}, 
\begin{equation*}
N_h\leq E_{h}(\hat{u})\leq \lim\limits_{k\to\infty}E_{h}(u_k)=N_h,   
\end{equation*}
which implies that \(u_k\to \hat{u}\) in \(H\). Then, \(\hat{u}\) is a minimizer of \(N_h\) and further a weak solution of Eq.\eqref{eq3.1}. The arguments in the proof of Lemma \ref{lm2.1} shows that \(\hat{u}\) is also a pointwise solution of Eq.\eqref{eq3.1}. Then by Eq.\eqref{eq3.1}, we have \(\|\hat{u}\|_H \neq 0\), since by \((h_1)\), \(\|h\|_2 \neq 0\), and therefore,
\begin{equation*}
\int_{V}h(x)\hat{u}d\mu=\int_{V}|\Delta\hat{u}|^2d\mu+\int_{V}(a(x)+1)|\hat{u}|^2d\mu=\| \hat{u}\|_{H}^2>0.   
\end{equation*}

\end{proof}

\begin{lemma}\label{lm3.2}
Let \(r_\varepsilon\) be given in Lemma \ref{lm3.1}, then there exists \(u'\in H\) independent of \(\varepsilon\), such that \(J_\varepsilon(u')<J_0(u')<0\) and \(\|u'\|_{H}>1>r_\varepsilon\).
\end{lemma}

\begin{proof}
Let \(\hat{u}\in H\) be given in Lemma \ref{lm-tech}. Since \(\int_{V}h(x)\hat{u}(x)d\mu>0\), then \(\hat{u}\) is not identically zero.  Denote \(A_t:=\{x\in V: t |\hat{u}(x)|>1\}\), then there exists $t_0>1$ such that 
\(meas(A_{t_0})>0\), where \(meas(A)\) denotes the Lebesgue measure of \(A\). Thus, for any \(t>t_0\), by \eqref{1.5} we deduce
\begin{equation}\label{3.1.1}
\begin{aligned}
J_\varepsilon(t\hat{u})&=\frac{1}{2}t^2\| \hat{u} \|_H^2-\int_{V}F(t\hat{u}(x))d\mu-t\int_{V}\varepsilon h(x)\hat{u}(x)d\mu+\frac{1}{2}t^2\|\nabla \hat{u}\|_{2}^2\\ 
&\leq\frac{1}{2}t^2\| \hat{u} \|_H^2-\beta\int_{A_{t_0}}t^\theta| \hat{u} |^\theta d\mu+\frac{1}{2}t^2\|\nabla \hat{u}\|_{2}^2\\ 
&=t^2\frac{1}{2}\{\| \hat{u} \|_H^2+\|\nabla \hat{u}\|_{2}^2\}-t^\theta\beta \int_{A_{t_0}}| \hat{u} |^\theta d\mu.
\end{aligned}
\end{equation}
In view of \(\theta>2\), we choose $\hat{t}>t_0$ sufficiently large, such that
\begin{equation}\label{3.1.2}
\hat{t}^2\frac{1}{2}\{\| \hat{u} \|_H^2+\|\nabla \hat{u}\|_{2}^2\}- \hat{t}^\theta\beta\int_{A_{t_0}}| \hat{u} |^\theta d\mu<0,\ \mbox{and}\ \|\hat{t}\hat{u}\|_H>1.
\end{equation} 
Note that \eqref{3.1.1} holds also as \(\varepsilon=0\). Thus, in view of \(\int_{V}\varepsilon h(x)\hat{u}(x)d\mu>0\), by letting $u'=\hat{t}\hat{u}$, we conclude that \(J_\varepsilon(u')<J_0(u')<0\) with \(\|u'\|_H>1\). Therefore, the proof is completed.
\end{proof}

\begin{lemma}\label{lm3.3}
The functional \(J_\varepsilon(u)\) satisfies the \((PS)_c\) condition. Namely, for a Palais-Smale sequence \(\{u_n\}\subset H\) of \(J_\varepsilon(u)\) satisfying: \(J_\varepsilon(u_n)\to c\) and \(J'_\varepsilon(u_n)\to0\), then up to a subsequence, there exists $u\in H$, such that \(u_n\to u\) in \(H\).
\end{lemma}

\begin{proof}
We first prove that \(\{u_n\}\) is bounded in \(H\). Let \(\theta>2\) be given in \eqref{1.4}, and observe that
\begin{align*}
 J_\varepsilon(u_n)-\frac{1}{\theta}(J'_\varepsilon(u_n),u_n)=\frac{\theta-2}{2\theta}\| u_n \|_{H}^2-\frac{1}{\theta}\int_{V}\Big(\theta F(u_n)- f(u_n)u_n\Big)d\mu\\
 -\frac{\theta-1}{\theta}\int_{V}\varepsilon h(x)u_nd\mu+\frac{\theta-2}{2\theta}\|\nabla  u_n\|_{2}^2.
\end{align*}
According to \eqref{1.4} and \(J_\varepsilon(u_n)\to c,J'_\varepsilon(u_n)\to 0\), we have
\begin{align*}
c+o_n(1)\| u_n \|_{H}=\frac{\theta-2}{2\theta}\| u_n \|_{H}^2-\frac{1}{\theta}\int_{V}\Big(\theta F(u_n)- f(u_n)u_n\Big)d\mu\\
 -\frac{\theta-1}{\theta}\int_{V}\varepsilon h(x)u_nd\mu+\frac{\theta-2}{2\theta}\|\nabla  u_n\|_{2}^2.
\end{align*}
By the H\"older inequality, we have:
\begin{equation}\label{3.6}
c+o_n(1)\| u_n \|_{H}+(\frac{\theta-1}{\theta})\varepsilon(a_0+1)^{-\frac{1}{2}}\|h(x)\|_{2}\|u_n\|_{H}\geq\frac{\theta-2}{2\theta}\| u_n \|_{H}^2,
\end{equation}
from which we deduce that \(\{u_n\}\) is bounded in \(H\). Thus by Lemma \ref{lm2.1} or Lemma \ref{lm2.2}, up to a subsequence, there exists \(u\in H\), such that
\begin{align}\label{3.1}
\begin{split}
u_n &\rightharpoonup u\ \mbox{in}\ H, \\
u_n &\to u\   \mbox{in}\ L^q(V),\ \forall q\in [2, +\infty].
\end{split}
\end{align}

Now we prove that \(u_n\to u\) in \(H\). Indeed, note that 
\begin{align*}
0\leftarrow(J'_\varepsilon(u_n),u_n-u)&=(u_n,u_n-u)_H-\int_{V}f( u_n)(u_n-u)d\mu\\
&-\int_{V}\varepsilon h(x) (u_n-u)d\mu+\int_{V}\Gamma(u_n,u_n-u)d\mu.
\end{align*}
From \eqref{1.3}, \eqref{3.1} and \((h_1)\), we have \(\lim\limits_{n \to \infty}(u_n,u_n-u)_{H}=\lim\limits_{n \to \infty}\| u_n \|_{H}^2-\| u \|_{H}^2\), and
\begin{align*}
\int_{V}f( u_n)(u_n-u)d\mu&\leq\int_{V}(\alpha u_n+C(\alpha)| u_n |^p )(u_n-u)d\mu\\
&\leq \alpha\| u_n \|_{2}\| u_n-u \|_{2}+C(\alpha)\|u_n-u\|_r\|u_n^p\|_{r'}\to 0,\\
\int_{V}\varepsilon h(x) (u_n-u)d\mu&\leq\varepsilon\|h(x)\|_{2}\|u_n-u\|_{2}\to 0,\\
\int_{V}\Gamma(u_n,u_n-u)d\mu&=\int_{V}\nabla u_n\nabla(u_n-u)d\mu\leq\|\Delta u_n\|_2\|u_n-u\|_2\to 0,
\end{align*}
where in the second inequality we have used the H\"older inequality with \(\dfrac{1}{r}+\dfrac{1}{r'}=1\) and \(r, r'>1\).

Hence, we conclude that 
\begin{equation*}
\lim\limits_{n \to \infty}\| u_n \|_{H}^2=\| u \|_{H}^2.
\end{equation*}
In view of its weak convergence, then \(u_n \to u\) in \(H\).
\end{proof}

For given \(\varepsilon>0\) small enough, let \(r_\varepsilon\) be given by Lemma \ref{lm3.1}, \(u'\) by Lemma \ref{lm3.2}, and define 
\begin{equation}\label{mp}
c_\varepsilon:=\inf\limits_{\gamma\in \mathcal{L}}\max\limits_{0\leq t\leq 1}J_\varepsilon(\gamma(t)),
\end{equation}
where the set of paths is given as:
\begin{equation}\label{Gamma}
\mathcal{L}:=\{\gamma\in C([0,1],H):\gamma(0)=0,\gamma(1)=u'\},
\end{equation}

Then, from Lemma \ref{lm3.1}, Lemma \ref{lm3.2} and Lemma \ref{lm3.3}, we deduce immediately the existence of a first solution.
\begin{lemma}\label{lm3.3.3}
For given \(\varepsilon>0\) small enough, there exists a critical point \(u_{\varepsilon}^{(1)}\in H\) of \(J_\varepsilon(u)\), which is of mountain-pass type and satisfies
\[J_\varepsilon(u_{\varepsilon}^{(1)})=c_\varepsilon>0.\]
\end{lemma}

Next, we show the existence of a second solution whose energy is negative and being of minimizer type. Let \(r_\varepsilon\) be given in Lemma \ref{lm3.1}, and define 
\begin{equation}\label{localmini}
m_\varepsilon:=\inf\limits_{\| u \|_{H}\leq \frac{1}{2}r_\varepsilon}J_\varepsilon(u).
\end{equation} 
Then, we can prove that 
\begin{lemma}\label{lm3.4}
For given \(\varepsilon>0\) small enough, there exists \(u_{\varepsilon}^{(2)}\in H\) with \(\|u_{\varepsilon}^{(2)}\|_H<\frac{1}{2}r_\varepsilon\), such that
\begin{equation*}
J_\varepsilon(u_{\varepsilon}^{(2)})=m_\varepsilon<0.
\end{equation*}
In particular, \(u_{\varepsilon}^{(2)}\) is a non-trivial solution of Eq. \eqref{eq1.1}.
\end{lemma}

\begin{proof}
Let \(\hat{u}\in H\) be given in Lemma \ref{lm-tech}. From \eqref{3.1.1} and Lemma \ref{lm-tech}, we have
\[J_\varepsilon(t\hat{u})\leq t^2\frac{1}{2}\{\| \hat{u} \|_H^2+\|\nabla \hat{u}\|_{2}^2\}-t\int_{V}\varepsilon h(x)\hat{u}(x)d\mu<0,\  \text{ and }\ \|t\hat{u}\|_H\leq \frac{1}{2}r_\varepsilon,\]
as \(t>0\) being chose sufficiently small. This gives that \(m_\varepsilon<0\).

Now, let \(\{u_k\}\subset H\) be a sequence such that 
\[\|u_k\|_H\leq\frac{1}{2} r_\varepsilon \ \ \mbox{and}\ \ \lim_{k\to \infty}J(u_k)= m_\varepsilon.\]
Then by Lemma \ref{lm2.1} or Lemma \ref{lm2.2} and their proofs, up to a subsequence, there exists \(u_{\varepsilon}^{(2)}\in H\) with \(\|u_{\varepsilon}^{(2)}\|_H\leq \frac{1}{2}r_\varepsilon\), such that
\begin{align*}
 u_k &\rightharpoonup u_{\varepsilon}^{(2)}, \ \text{in }H,\\
  u_k &\to u_{\varepsilon}^{(2)}, \  \text{in }L^q(V),\ \ \forall q\in [2, +\infty];\\
  u_k(x)&\to u_{\varepsilon}^{(2)}(x), \quad \forall x \in V.
\end{align*}
Moreover,
\begin{align*}
\lim\limits_{k \to \infty}\int_{V}\varepsilon h(x)(u_k-u_{\varepsilon}^{(2)})d\mu &\leq\lim\limits_{k \to \infty}\varepsilon\|h\|_2\|u_k-u_{\varepsilon}^{(2)}\|_2\to 0,\\
\lim\limits_{k \to \infty}\int_{V}\Big | |\nabla u_k|^2-|\nabla u_{\varepsilon}^{(2)}|^2 \Big |d\mu&=\lim\limits_{k \to \infty}\int_{V}|\nabla(u_k-u_{\varepsilon}^{(2)})||\nabla(u_k+u_{\varepsilon}^{(2)})|d\mu\\
&\leq\lim\limits_{k \to \infty}\|\Delta(u_k+u_{\varepsilon}^{(2)})\|_2\|u_k-u_{\varepsilon}^{(2)}\|_2\to 0,\\
\|u_{\varepsilon}^{(2)}\|_{H}&\leq\lim\limits_{k \to \infty}\inf\|u_k\|_{H}\leq\frac{1}{2}r_\varepsilon.
\end{align*}
By the Lebesgue dominated convergence theorem, we obtain:
\begin{equation*}
\lim\limits_{k \to \infty}\int_{V}F(u_k)d\mu=\int_{V}F(u_{\varepsilon}^{(2)})d\mu,
\end{equation*}
Therefore
\begin{equation*}
m_\varepsilon\leq J(u_{\varepsilon}^{(2)})\leq \lim\limits_{k \to \infty}\inf J(u_k)=m_\varepsilon,   
\end{equation*}
which implies that \(u_{\varepsilon}^{(2)}\) is a minimizer of \(m_\varepsilon\).  In particular, since \(J_{\varepsilon}(u_{\varepsilon}^{(2)})=m_\varepsilon<0\), then from Lemma \ref{lm3.1} we know that \(\|u_{\varepsilon}^{(2)}\|_{H}<\frac{1}{2}r_\varepsilon\) and \(u_{\varepsilon}^{(2)}\neq 0\).

Now, we show that \(u_{\varepsilon}^{(2)}\) is indeed a non-trivial weak solution of Eq.\eqref{eq1.1}. For any \(\varphi\in C_c(V)\), let \[\psi(t):=J_\varepsilon(u_{\varepsilon}^{(2)}+t\varphi), \ t\in \R.\] 
Then \(\psi\) is a smooth function on \(\R \). Since \(\|u_{\varepsilon}^{(2)}\|_H<\frac{1}{2}r_\varepsilon\), then we can choose \(\eta>0\) small enough, such that \(\|u_{\varepsilon}^{(2)}+t\varphi\|_H\leq \frac{1}{2}r_\varepsilon\) if \(|t|<\eta \). Using the fact that \(u_{\varepsilon}^{(2)}\) is the minimizer of \(m_\varepsilon\), we have \(\psi^{\prime}(0)=0\), which gives that
\begin{equation*}
\int_{V}[\Delta u_{\varepsilon}^{(2)}\Delta\varphi+\nabla u_{\varepsilon}^{(2)}\nabla\varphi+(a(x)+1)u_{\varepsilon}^{(2)}\varphi-f(u_{\varepsilon}^{(2)})\varphi-\varepsilon h(x)\varphi]d\mu=0.
\end{equation*}
This shows that \(u_{\varepsilon}^{(2)}\) is a non-trivial weak solution of Eq.\eqref{eq1.1}. Thus the proof is completed.
\end{proof}

Now, we are ready to prove Theorem \ref{th1.1}.
\begin{proof}[\textbf{Proof of Theorem \ref{th1.1}}]
In view of Lemma \ref{lm3.3.3} and Lemma \ref{lm3.4}, we have established the existence of two non-trivial solutions of Eq.\eqref{eq1.1}, where one is of mountain pass type with positive energy and another is of minimizer type with negative energy. 

%
\end{proof}

In the following, we prove Theorem \ref{th1.2}. First, we show that
\begin{lemma}\label{lm3.5}
For given \(r\in (0,1)\) being small enough, there exists constant \(L_r>0\), such that \(J_0(u)\geq L_r\) for all  \(\|u\|_H=r\).
\end{lemma}

\begin{proof}
Let \(\alpha=\frac{a_0+1}{2}\) in \eqref{1.3}, then by Remark \ref{rek2.1} and \eqref{1.3}, we have
\begin{equation*}
\int_{V}F(u)d\mu\leq\frac{1}{4}\| u \|_{H}^2+\frac{C_{p+1}}{p+1}C(\alpha)\| u \|_{H}^{p+1}, \ \forall  u\in H.
\end{equation*}
Then
\begin{equation*}
\begin{aligned}
J_0(u)&\geq\frac{1}{2}\| u \|_{H}^2+\frac{1}{2}\|\nabla u\|_{2}^2-\frac{1}{4}\| u \|_{H}^2-\frac{C_{p+1}}{p+1}C(\alpha)\| u \|_{H}^{p+1}\\&\geq\frac{1}{4}\| u \|_{H}^2-\frac{C_{p+1}}{p+1}C(\alpha)\| u \|_{H}^{p+1}.
\end{aligned}
\end{equation*}
Let \(L_r=\frac{1}{4}r^2-\frac{C_{p+1}}{p+1}C(\alpha)r^{p+1}\), then by letting \(r>0\) small enough with also \(r<1\), such that \(L_r>0\). Then it follows that \(J_0(u) \geq L_r\) for all \(\|u\|_H = r\).
\end{proof}

\begin{lemma}\label{lm3.6}
Let \(u'\in H\) be given in Lemma \ref{lm3.2}, then \(J_0(u')<0\) and \(\|u'\|_{H}>r\).
\end{lemma}
\begin{proof}
This lemma follows immediately from Lemma \ref{lm3.2}.
\end{proof}

Let \(\mathcal{L}\) be given in \eqref{Gamma}, then we define
\begin{equation}\label{mp_0}
c_0:=\inf\limits_{\gamma\in\mathcal{L}}\max\limits_{0\leq t\leq 1}J_0(\gamma(t)).
\end{equation}
Then, we show that \(c_0>0\) is a critical value of the functional \(J_0\). 
\begin{proof}[\textbf{Proof of Theorem \ref{th1.2}}]
Observe that in the proof of Lemma \ref{lm3.3} \(\varepsilon>0\) does not enter essentially the analysis, then using the same arguments with \(\varepsilon=0\), we can verify that the functional \(J_0(u)\) also satisfies the \((PS)_c\) condition. Then, in view of Lemma \ref{lm3.5} and Lemma \ref{lm3.6}, there exists \(u_0\in H\) which is a critical point of \(J_0\) with \(J_0(u_0)=c_0>0\). 
\end{proof}

\section{Proofs of Theorem \ref{th1.3}}
In this section, we prove the asymptotic behaviors of the local minimum solution \(u_{\varepsilon}^{(2)}\) and the mountain pass solution \(u_{\varepsilon}^{(1)}\) when \(\varepsilon\to 0\), where these two solutions have been obtained by Theorem \ref{th1.1} and Theorem \ref{th1.2}.

\begin{lemma}\label{lm4.1}
As \(\varepsilon\to 0^+\), the local minimum solution \(u_{\varepsilon}^{(2)}\) of Eq.\eqref{eq1.1} satisfies \(u_{\varepsilon}^{(2)}\to 0\) in \(H\), and \(J_\varepsilon(u_{\varepsilon}^{(2)})\to 0\).
\end{lemma}
\begin{proof} 
By Lemma \ref{lm3.4}, we have \(\|u_{\varepsilon}^{(2)}\|_H<\frac{1}{2}r_\varepsilon=\frac{1}{2}\sqrt{\varepsilon}\). Then,
\begin{equation*}
0\leq\lim\limits_{\varepsilon\to 0^+}\|u_{\varepsilon}^{(2)}\|_H\leq\lim\limits_{\varepsilon\to 0^+}\frac{1}{2}\sqrt{\varepsilon}= 0,
\end{equation*}
which gives
\begin{equation}
\lim\limits_{\varepsilon\to 0^+}\|u_{\varepsilon}^{(2)}-0\|_H=0.   
\end{equation}
From \eqref{estima} and \(\|u_{\varepsilon}^{(2)}\|_H\leq \frac{1}{2}\sqrt{\varepsilon}\), we have
\begin{equation*}
\begin{aligned}
J_\varepsilon(u_{\varepsilon}^{(2)})&\geq\frac{1}{4}\|u_{\varepsilon}^{(2)}\|_{H}^2-\frac{C_{p+1}}{p+1}C(\alpha)\|u_{\varepsilon}^{(2)}\|_{H}^{p+1}-\varepsilon (a_0+1)^{-\frac{1}{2}}\|h(x)\|_{2}\|u_{\varepsilon}^{(2)}\|_{H}\\&\geq-\frac{C_{p+1}}{p+1}C(\alpha)(\frac{\sqrt{\varepsilon}}{2})^{p+1}-\frac{1}{2}(a_0+1)^{-\frac{1}{2}}\varepsilon^{\frac{3}{2}}\|h\|_2.
\end{aligned}
\end{equation*}
Take the limit of \(\varepsilon \to 0^+\), then
\begin{equation}\label{4.2}
\liminf\limits_{\varepsilon\to 0^+}J_\varepsilon(u_{\varepsilon}^{(2)})\geq0.   
\end{equation}
On the other hand, from Lemma \ref{lm3.4}, we know that
\begin{equation}\label{4.3}
J_\varepsilon(u_{\varepsilon}^{(2)})=m_\varepsilon<0,   
\end{equation}
Combining \eqref{4.2} and \eqref{4.3}, we conclude that
\begin{equation*}
\lim\limits_{\varepsilon\to 0^+}J_\varepsilon(u_{\varepsilon}^{(2)})=0.
\end{equation*}
This completes the proof.
\end{proof}

\begin{lemma}\label{lm4.2}
As \(\varepsilon\to 0^+\), the mountain pass solution \(u_{\varepsilon}^{(1)}\) of Eq.\eqref{eq1.1} satisfies that 
\[u_{\varepsilon}^{(1)}\to u_1,\ \mbox{ in } H, \] 
where \(u_1\in H\) is a mountain pass solution of \eqref{eq1.2} at the energy level \(J_0(u_1)=c_0\). In particular, 
\[c_\varepsilon=J_{\varepsilon}(u_{\varepsilon}^{(1)}) \to J_0(u_1),\] 
where \(c_{\varepsilon}, c_0\) are given by \eqref{mp} and \eqref{mp_0} respectively. 
\end{lemma}

\begin{proof} First, we prove that 
\begin{eqnarray}\label{4.4}
c_\varepsilon\to c_0, \ \mbox{as}\ \varepsilon\to 0^+.
\end{eqnarray}
To this aim, on one hand, we prove 
\begin{eqnarray}\label{upper}
\limsup_{\varepsilon \to 0^+} c_\varepsilon \leq c_0.
\end{eqnarray}
In fact, by definition of \(c_0\), for any \(\delta > 0\), there exists a path \(\gamma_\delta \in \mathcal{L}\) such that
\[
\max_{0\leq t\leq 1} J_0(\gamma_\delta(t)) \leq c_0 + \delta.
\]
For this fixed path, compute the difference:
\[
|J_\varepsilon(\gamma_\delta(t)) - J_0(\gamma_\delta(t))| = \left| \varepsilon \int_V h(x) \gamma_\delta(t) d\mu \right| \leq \varepsilon \|h\|_2 \|\gamma_\delta(t)\|_2.
\]
Let \(M_\delta := \|h\|_2 \cdot \max_{t \in [0,1]} \|\gamma_\delta(t)\|_2\) (independent of \(\varepsilon\)), then
\[
c_{\varepsilon}\leq \max_{0\leq t\leq 1} J_\varepsilon(\gamma_\delta(t)) \leq c_0 + \delta + \varepsilon M_\delta.
\]
Taking \(\limsup\limits_{\varepsilon \to 0^+}\) and using the arbitrariness of \(\delta\), we get immediately \eqref{upper}.

On the other hand, we prove the inverse inequality
\begin{eqnarray}\label{lower}
\liminf_{\varepsilon \to 0^+} c_\varepsilon \geq c_0.
\end{eqnarray}
For this purpose, we write
\[
A(u):=\int_V\big(|\Delta u|^2+|\nabla u|^2+(a(x)+1)u^2\big)d\mu,
\quad \ell(u):=\int_Vh(x)u(x)d\mu.
\]
Thus $J_\varepsilon(u)=\frac12A(u)-\int_VF(u)d\mu-\varepsilon\ell(u)$. By \((A_1)\), with $K:=(a_0+1)^{-\frac{1}{2}}\|h\|_2$,
\begin{equation}\label{eq:ell}
|\ell(u)|\le K\|u\|_H\le K\sqrt{A(u)},\quad \forall u\in H.
\end{equation}

Note that \((f3)\) implies $sf(s)\ge \theta F(s)$ for every $s\in \R$, where \(\theta>2\).
Thus
\[
\frac{d}{d\lambda}\left(\lambda^{-\theta}F(\lambda s)\right)
=\lambda^{-\theta-1}\bigl[\lambda sf(\lambda s)-\theta F(\lambda s)\bigr]\ge0,
\]
and hence
\begin{equation}\label{eq:scaling}
F(\lambda s)\ge\lambda^{\theta}F(s),\quad(\lambda\ge1,\ s\in\mathbb R).
\end{equation}
Fix $\lambda>1$. By \eqref{eq:ell}, \eqref{eq:scaling}, and the Young inequality,
\begin{align}\label{eq:comparison}
J_0(\lambda u)
&\le\lambda^{\theta}J_\varepsilon(u)
-\frac{\lambda^{\theta}-\lambda^2}{2}A(u)
+\varepsilon\lambda^{\theta}\ell(u)\notag\\
&\le\lambda^{\theta}J_\varepsilon(u)
+\frac{\varepsilon^2\lambda^{2\theta}K^2}
{2(\lambda^{\theta}-\lambda^2)},
\quad(u\in H).
\end{align}
Let  $u'$ be given by Lemma \ref{lm3.2}, which is independent of $\varepsilon$, then $J_0(u')=\lim_{\varepsilon\to0^+}J_\varepsilon(u')\le0$.
For every $s\ge1$, \eqref{eq:scaling} gives
\begin{equation}\label{eq:connector}
J_0(su')\le s^{\theta}J_0(u')-\frac{s^{\theta}-s^2}{2}A(u')\le0.
\end{equation}
For any $\gamma\in \mathcal{L}$, define
\[
\widetilde\gamma(t):=
\begin{cases}
\lambda\gamma(2t),&0\le t\le\frac12,\\
[\lambda-2(\lambda-1)(t-\frac12)]u',&\frac12\le t\le1.
\end{cases}
\]
Then $\widetilde\gamma\in \mathcal{L}$. Since $J_\varepsilon(\gamma(0))=0$,
the maximum of $J_\varepsilon$ along $\gamma$ is nonnegative.
Using \eqref{eq:comparison} and \eqref{eq:connector}, we obtain
\[
c_0\le\max_{0\leq t\leq 1}J_0(\widetilde\gamma(t))
\le\lambda^{\theta}\max_{0\leq t\leq 1}J_\varepsilon(\gamma(t))+\frac{\varepsilon^2\lambda^{2\theta}K^2}
{2(\lambda^{\theta}-\lambda^2)}.
\]
Taking the infimum over $\gamma\in \mathcal{L}$ yields
\[
c_\varepsilon\ge\lambda^{-\theta}c_0
-\frac{\varepsilon^2K^2}{2(1-\lambda^{2-\theta})}.
\]
First let $\varepsilon\to0^+$ with $\lambda>1$ fixed, and then let
$\lambda\downarrow1$. We obtain \eqref{lower}. Combining \eqref{upper} and \eqref{lower}, then we have verified \eqref{4.4}. \\

Next, we prove that \(u_{\varepsilon}^{(1)}\) is uniformly bounded in \(H\).  Indeed, from Lemma \ref{lm3.3.3}, we see that \(J_\varepsilon(u_{\varepsilon}^{(1)})=c_\varepsilon\) and \((J'_\varepsilon(u_{\varepsilon}^{(1)}),u_{\varepsilon}^{(1)})=0\).
Using these two facts, and the arguments in the proof in Lemma \ref{lm3.3}, we obtain as \eqref{3.6} that
\begin{equation*}
c_\varepsilon+(1-\frac{1}{\theta})(a_0+1)^{-\frac{1}{2}}\varepsilon\|h(x)\|_{2}\|u_{\varepsilon}^{(1)}\|_{H}\geq\frac{\theta-2}{2\theta}\| u_{\varepsilon}^{(1)} \|_{H}^2, \ \theta>2.
\end{equation*}
In view of \(c_\varepsilon\to c_0\), then \(\|u_{\varepsilon}^{(1)}\|_H\) is bounded uniformly. Since \(H\) is reflexive and Lemma \ref{lm2.1} or Lemma \ref{lm2.2}, then there exists \(u_1\in H\).
\begin{equation}\label{4.5}
u_{\varepsilon}^{(1)} \rightharpoonup u_1 \quad \text{in }H,\qquad u_{\varepsilon}^{(1)}\to u_1\quad \text{in }L^q,\quad (2 \leq q \leq +\infty  ).   
\end{equation}
Thus, by \(J'_\varepsilon(u_{\varepsilon}^{(1)})=0\), we have
\begin{eqnarray*}
0=(J'_\varepsilon(u_{\varepsilon}^{(1)}),u_{\varepsilon}^{(1)}-u_1) =(u_{\varepsilon}^{(1)},u_{\varepsilon}^{(1)}-u_1)_H-\int_{V}f( u_{\varepsilon}^{(1)})(u_{\varepsilon}^{(1)}-u_1)d\mu\\
-\int_{V}\varepsilon h(x) (u_{\varepsilon}^{(1)}-u_1)d\mu+\int_{V}\Gamma(u_{\varepsilon}^{(1)},u_{\varepsilon}^{(1)}-u_1)d\mu.
\end{eqnarray*}
As the proof in Lemma \ref{lm3.3}, we get from the above that \((u_{\varepsilon}^{(1)},u_{\varepsilon}^{(1)}-u_1)_H \to 0 \) as \(\varepsilon \to 0^+\). This leads to 
\(\|u_{\varepsilon}^{(1)}\|_H \to \|u_1\|_H\), and therefore 
\begin{equation}\label{4.5.5}
u_{\varepsilon}^{(1)} \to u_1\quad \mbox{in}\ H.
\end{equation}

Since \eqref{4.5.5}, we have \(J_0(u_{\varepsilon}^{(1)}) \to J_0(u_1)\). Moreover, note that 
\[
|J_\varepsilon(u_{\varepsilon}^{(1)}) - J_0(u_{\varepsilon}^{(1)})| = \left| \varepsilon \int_V h(x) u_{\varepsilon}^{(1)} d\mu \right| \leq \varepsilon \|h\|_2 \|u_{\varepsilon}^{(1)}\|_2\leq \varepsilon (a_0+1)^{-\frac{1}{2}}\|h\|_2 \|u_{\varepsilon}^{(1)}\|_H.
\]
Then, in view of the fact that \(u_{\varepsilon}^{(1)}\) is uniformly bounded in \(H\), we obtain  
\[c_{\varepsilon}=J_\varepsilon(u_{\varepsilon}^{(1)}) = J_0(u_{\varepsilon}^{(1)})+o_{\varepsilon}(1)=J_0(u_1)+o_{\varepsilon}(1).\] 
In particular, by \eqref{4.4}, \(J_0(u_1)=c_0\).\\

Finally, we prove that \(u_1\) is a weak solution to Eq.\eqref{eq1.2}. From Lemma \ref{lm3.3.3}, for any \(\varphi\),
\begin{equation}\label{4.6}
\int_{V}[\Delta u_{\varepsilon}^{(1)}\Delta\varphi+\nabla u_{\varepsilon}^{(1)}\nabla\varphi+(a(x)+1)u_{\varepsilon}^{(1)}\varphi-f(u_{\varepsilon}^{(1)})\varphi-\varepsilon h(x)\varphi]d\mu=0. 
\end{equation}
Take \(\varepsilon\to 0^+\). By \eqref{4.5}, we have
\begin{equation*}
\int_{V}\Delta u_{\varepsilon}^{(1)}\Delta\varphi\to\int_{V}\Delta u_1\Delta\varphi,
\end{equation*}
\begin{equation*}
\int_{V}(a(x)+1)u_{\varepsilon}^{(1)}\varphi\to\int_{V}(a(x)+1)u_1\varphi,
\end{equation*}
\begin{equation*}
\int_{V}\nabla (u_{\varepsilon}^{(1)}-u_1)\nabla\varphi d\mu\leq\|\Delta\varphi\|_2\|u_{\varepsilon}^{(1)}-u_1\|_2\to 0,
\end{equation*}
By the Lebesgue dominated convergence theorem, we have:
\begin{equation*}
\int_{V}f(u_{\varepsilon}^{(1)})\varphi d\mu\to\int_{V}f(u_1)\varphi d\mu,
\end{equation*}
In summary,
\begin{equation}
\int_{V}[\Delta u_1\Delta\varphi+\nabla u_1\nabla\varphi+(a(x)+1)u_1\varphi-f(u_1)\varphi]d\mu=0.   
\end{equation}
This completes the proof.
\end{proof}

\begin{proof}[\textbf{Proof of Theorem \ref{th1.3}}]
Clearly, Theorem \ref{th1.3} follows directly from Lemma \ref{lm4.1} and Lemma \ref{lm4.2}.
\end{proof}

%

\section*{Declarations}

\noindent\textbf{Funding} \\
Tingjian Luo is supported by the National Natural Science Foundation of China (11501137), and the Guangdong Basic and Applied Basic Research Foundation (2025A1515012282).

\medskip

\noindent\textbf{Author Contributions} \\
\textbf{Tingjian Luo:} Initiated the main ideas, reviewed and edited the final version of the manuscript.
\textbf{Yongqiang Wei:} drafted the initial manuscript and checked the proofs.
\textbf{All authors} approved the final submission.

\medskip

\noindent\textbf{Data Availability Statement} \\
No new data were created or analyzed in this study. Data sharing is not applicable to this article.

\medskip

\noindent\textbf{Conflict of Interest} \\
The authors declare no conflict of interest.


\begin{thebibliography}{99}
\bibitem{AR} Ambrosetti A., Rabinowitz P.H., Dual variational methods in critical point theory and applications, \textit{J. Func. Anal.} {\bf 14} (1973) 349–381.

\bibitem{BL} Bose N. K., Liang P., Neural network fundamentals with graphs, algorithms, and applications, McGraw-Hill, Inc. 1996.

\bibitem{BCGJ1} Bonheure, D., Casteras, J.-B., Gou, T., Jeanjean, L., Normalized solutions to the mixed dispersion nonlinear Schr\"odinger equation in the mass critical and supercritical regime. \textit{Trans. Am. Math. Soc.} {\bf 372} (2019), 2167-2212.

\bibitem{BCGJ2} Bonheure,D.,Casteras, J.-B., Gou, T., Jeanjean, L., Strong instability of ground states to a fourth order Schr\"odinger equation. \textit{Int. Math. Res. Not.} {\bf 17} (2019), 5299–5315.
\bibitem{BCM} Bonheure, D., Casteras, J.-B., Mandel, R., On a fourth-order nonlinear Helmholtz equation. \textit{J. Lond. Math. Soc.} {\bf 99} (2019), 831-852.

\bibitem{BCSN} Bonheure, D., Casteras, J.-B., dos Santos, E.-M., Nascimento, R., Orbitally stable standing waves of a mixed dispersion nonlinear Schr\"odinger equation. \textit{SIAM J. Math. Anal.} {\bf 50} (2018), 5027-5071.

\bibitem{BL} Boulenger, T., Lenzmann, E., Blowup for biharmonic NLS. \textit{Ann. Sci. Éc. Norm. Supér.}  {\bf 50} (2017), 503-544.

\bibitem{BN} Bonheure, D., Nascimento, R., Waveguide solutions for a nonlinear Schr\"odinger equation with mixed dispersion. In: Contributions to Nonlinear Elliptic Equations and Systems, Volume 86 of Progr. Nonlinear Differential Equations Appl., Birkhäuser/Springer, Cham (2015) pp. 31-53. 

\bibitem{CF} Casteras, J. B., Földes, J., Existence of traveling waves for a fourth order Schr\"odinger equation with mixed dispersion in the Helmholtz regime. \textit{Math. Ann.} {\bf 41} (2026). https://doi.org/10.1007/s00208-026-03399-4.

\bibitem{CM} Camilli F., Marchi C.,  A note on kazdan–warner equation on networks, \textit{Adv. Calc. Var.}, {\bf 15} (2022), 693-704. 
 
\bibitem{CWY} Chang X., Wang R., Yan D., Ground states for logarithmic Schr\"odinger equations on locally finite graphs,  \textit{J. Geom. Anal.}, {\bf 33} (2023), 26 pp.
    
\bibitem{FJMM} Fernández A. J., Jeanjean, L., Mandel R., Mariş M., Non-homogeneous Gagliardo-Nirenberg inequalities in $\bold{R}^N$ and application to a biharmonic non-linear Schr\"odinger equation, \textit{J. Differential Equations} {\bf 330} (2022), 1-65.
        
\bibitem{LW} Lin Y., Wu Y., The existence and nonexistence of global solutions for a semilinear heat equation on graphs,  \textit{Calc. Var. Partial Differential Equations}, {\bf 56} (2017), 22 pp.
    
\bibitem{LZ} Li Y., Zhou Z., Asymptotic Behavior of Solutions to Nonlinear Schrödinger Equations with Perturbations on Graphs, \textit{Bull. Malays. Math. Sci. Soc.}, accepted, 2026.

\bibitem{G} Ge H., The $p$-th Kazdan-Warner equation on graphs, \textit{Commun. Contemp. Math.}, \textbf{22} (2020), no.~6, 1950052, 17 pp. 

\bibitem{GHZ} Guo Z., Huang X., Zhou F.,  Radial symmetry of entire solutions of a bi-harmonic equation with exponential nonlinearity,  \textit{J. Funct. Anal.}, \textbf{268} (2015), 1972–2004. 
    
\bibitem{HSZ} Han X., Shao M., Zhao L., Existence and convergence of solutions for nonlinear bi-harmonic equations on graphs, \textit{J. Differential Equations},  \textbf{268} (2020), 3936-3961.

\bibitem{H} Hou S., Multiple Solutions of a Nonlinear Biharmonic Equation on Graphs, \emph{Commun. Math. Stat.}, \textbf{11} (2023), 767-774. 

\bibitem{K} Karpman V. I., Stabilization of soliton instabilities by higher order dispersion: Fourth-order nonlinear Schr\"odinger-type equations, \emph{Phys. Lett. E}, \textbf{53} (1996), 1336-1339.
\bibitem{KS} Karpman V. I., Shagalov A. G., Stability of solitons described by nonlinear Schr\"odinger-type equations with higher-order dispersion, \emph{Phys. D.}, \textbf{144} (2000), 194-210.

\bibitem{LY} Liu Y., Yang Y.,  Topological degree for kazdan–warner equation in the negative case on finite graph,  \textit{Ann. Global Anal. Geom.}, {\bf 65} (2024), 20 pp.
    
\bibitem{S} Sani F., A bi-harmonic equation in $\Bbb R^4$ involving nonlinearities with critical exponential growth, \textit{Commun. Pure Appl. Anal.}, {\bf 12} (2013), 405-428. 

\bibitem{SW} Sun L., Wang L., Brouwer degree for Kazdan-Warner equations on a connected finite graph, \textit{Adv. Math.}, {\bf 404} (2022), 29 pp. 

\bibitem{WS} Wang Y., Shen Y.,  Multiple and sign-changing solutions for a class of semilinear bi-harmonic equation,  \textit{J. Differential Equations}, \textbf{246} (2009), 3109–3125.

\bibitem{YLY1} Grigor'yan A., Lin Y., Yang Y., Kazdan-Warner equation on graph, \textit{Calc. Var. Partial Differential Equations}, {\bf 55} (2016), 13 pp. 
    
\bibitem{YLY2} Grigor'yan A., Lin Y., Yang Y.,  Existence of positive solutions to some nonlinear equations on locally finite graphs,  \textit{Sci. China Math.}, \textbf{60} (2017), 1311–1324.

\bibitem{ZLJS} Zhang J., Lou Z., Ji Y., Shao W.,  Multiplicity of solutions of the bi-harmonic Schr\"odinger equation with critical growth, \textit{Z. Angew. Math. Phys.}, \textbf{69} (2018).

\bibitem{ZZ} Zhao L., Zhang N.,  Existence of solutions for a higher order kirchhoff type problem with exponential critical growth, \textit{Nonlinear Analysis}, \textbf{132} (2016), 214-226.
    
\end{thebibliography}
\end{document}